

\input amstex.tex

\documentstyle{amsppt}
\magnification=1200
\hsize=5.0in\vsize=7.0in
\hoffset=0.2in\voffset=0cm
\nonstopmode
\nologo
\TagsAsMath


\def\det{\operatorname{det}}

\def\rank{\operatorname{rank}}

\def\ker{\operatorname{Ker}}
\def\supp{\operatorname{supp}}
\def\sup{\operatorname{sup}}
\def\inf{\operatorname{inf}}

\def\p{\partial}\def\at#1{\vert\sb{\sb{#1}}}

\def\vect#1{{\bold #1}}
\def\R{{\Bbb R}} 
 
\def\N{{\Bbb N}} 
\def\Z{{\Bbb Z}}

\def\Abs#1{\left\vert#1\right\vert}
\def\abs#1{\vert#1\vert}
\def\Norm#1{\left\Vert #1 \right\Vert}
\def\norm#1{\Vert #1 \Vert}

\def\ssb#1{\!\,\sb{\!\sb{#1}}}

\def\csb#1{\sb{\hskip -1pt #1}}
\def\cpr{\sp{\hskip -0.5pt {\sssize \prime}}}

\def\const{\operatorname{const}}

\def\nut{\sb{\!o}}
\def\loc{\sb{\text{loc}}}
\def\comp{\sb{\text{comp}}}

\def\authorinfo{%
\author
Andrew Comech
\medskip
{\sl Department of Mathematics\\ 
University of Toronto\\
Toronto, Canada {\it M5S\,3G3}}
\endauthor
\email { comech${\sssize\@}$math.toronto.edu}\endemail
}%

\def\sec#1{#1}

\def\ye#1{$\overset{\ssize #1}\to{\phantom{.}}$}
\def\noref#1\endref{}
\def\hvar{\hslash}

\def\thvar{\vartheta}

\def\sing#1{{S}\sb{#1}}

\document

\topmatter

\title
Damping estimates
for
oscillatory integral operators
with finite type singularities
\endtitle

\authorinfo


\leftheadtext{Andrew Comech}
\rightheadtext{Damping estimates for oscillatory integral operators}

\abstract
We derive 
damping estimates
and 
asymptotics of $L^p$ operator norms
for oscillatory integral operators
with finite type singularities.
The methods are based on incorporating 
finite type conditions into $L^2$
almost orthogonality technique of Cotlar-Stein.
\endabstract

\endtopmatter 

\head \sec1. Introduction and results
\endhead

The oscillatory integral operators
have the form
$$
T\sb{\lambda}u(x)
=\underset{\R^n}\to\int
e^{i\lambda S(x,\vartheta)}
\psi(x,\vartheta)
u(\vartheta)\,d\vartheta,
\qquad
x,\,\vartheta\in\R^n,
\tag{\sec1.1}
$$
with $S\in C^\infty(\R^n\times\R^n)$
and $\psi\in C\comp^\infty(\R^n\times\R^n)$.
We denote $h(x,\vartheta)=\det S\sb{x\vartheta}$.
If $h(x,\vartheta)\ne 0$, then it is well-known 
that the $L^2$ operator norm of $T\sb{\lambda}$
decays as $\lambda^{-\frac{n}{2}}$ \cite{H\"o\ye{71}}.
The operators with non-empty 
critical variety
$$
\varSigma=\{(x,\vartheta)\,\,\vert\,\, h(x,\vartheta)=0\}
$$
attracted much attention
during last several years:
\cite{PaSo\ye{90}}, \cite{Pa\ye{91}},
\cite{PhSt\ye{91}} -- \cite{PhSt\ye{97}},
\cite{GrSe\ye{94}} -- \cite{GrSe\ye{97b}},
\cite{Cu\ye{97}}.
We recommend \cite{Ph\ye{95}}
as a survey on integral operators
associated to singular canonical relations.

The properties of $T\sb{\lambda}$ are being
characterized in terms of the projections
from the associated canonical relation
$
{\Cal C}
=\{(x,\,S\sb{x},\,\thvar,\,S\sb{\thvar})\}
\subset T\sp{\ast}\R^n\times T\sp{\ast}\R^n
$
onto the left and right factors.
We consider these projections as lifted onto
$\R^n\times\R^n$:
$$
\pi\ssb{L}:\;(x,\thvar)\mapsto (x,\,S\sb{x}),
\qquad
\pi\ssb{R}:\;(x,\thvar)\mapsto (\thvar,\,S\sb{\thvar}).
$$
These maps become singular 
on the critical variety $\varSigma$.

We already know \cite{Co\ye{98}}
that if 
one of the projections from the canonical relation
is a Whitney fold 
while the type of the other projection
is at most $k$ ($k=1$ for a Whitney fold),
then 
$\norm{T\sb{\lambda}}
\le\const\lambda^{-\frac{n}{2}+(4+\frac{2}{k})^{-1}}$.
This result was used to obtain the 
optimal regularity of Fourier integral
operators with one-sided Whitney folds.

In this paper, we approach a much more
complicated situation when both projections
from the canonical relation are of finite type.
We develop the machinery 
which yields the asymptotics of 
the norm of $(\sec1.1)$
with the integral kernel localized to the 
region where $h(x,\vartheta)\sim\hvar$,
$\hvar$ being some small real number.
We then derive 
the damping estimates 
on oscillatory integral operators:
we will prove that
if the symbol of the operator
vanishes as $\abs{h(x,\vartheta)}$
on the critical variety,
then the operator has such  properties 
as though it is associated to a canonical graph:
$\norm{T\sb{\lambda}}\le\const\lambda^{-\frac{n}{2}}$.
This result was previously proved 
for operators associated to 
two-sided Whitney folds \cite{MeT\ye{85}}
and for operators in $n=1$, with polynomial
phases \cite{PhSt\ye{94}}.
Damping 
for operators with one-sided Whitney folds 
follows from \cite{Co\ye{97}}.
A much more general situation 
(when no assumptions on the projections
$\pi\ssb{L}$, $\pi\ssb{R}$ are made)
is considered in \cite{SoSt\ye{86}}:
the damping occurs if the symbol vanishes
as $\abs{h(x,\vartheta)}^{5n/2}$.

We will exploit the concept of
the type of a map,
which we define as 
the highest order of vanishing
of the determinant of its Jacobi matrix 
in the ``critical'' direction
\cite{Co\ye{98}}.
Let $M$ and $N$ be two $C^\infty$ manifolds 
of the same dimension
and let $\pi:\;M\rightarrow N$
be a smooth map
with corank at most 1.
Assume that $\det d\pi$ vanishes simply 
on $\varSigma\subset M$.

\proclaim{Definition}
Let $\vect{V}$
be any smooth vector field 
which generates (locally) the kernel of $d\pi$:
$
\vect{V}\at{\varSigma}\in\ker d\pi,
$
$
\vect{V}\at{\varSigma}\ne 0.
$
The type of $\pi$ at a point
$p\nut\in\varSigma$
is defined to be  
the smallest integer $k$ such that
$
\vect{V}^{k} \det d\pi\at{p\nut}\ne 0.
$

The type of $\pi$ at $p\notin\varSigma$
is defined to be $0$.
\endproclaim
\noindent
An example of a map of type at most $k$
is a map which has a Morin 
$\sing{1\csb{k}}$-singularity
\cite{Mo\ye{65}}.
In particular,
the Whitney fold is of type at most $1$.

\subhead
Asymptotics of $L^2$ estimates
\endsubhead
Let us localize the integral kernel 
of $T\sb{\lambda}$
with the aid of a certain smooth function
$\beta$
to the region where $\det S\sb{x\vartheta}$
takes the values of size $\hvar$
(for simplicity, we assume that $\hvar>0$):
$$
T\sb{\lambda}\sp{\hvar}u(x)
=\underset{\R^n}\to\int
e^{i\lambda S(x,\vartheta)}
\psi(x,\vartheta)
\beta(\hvar^{-1} h(x,\vartheta))
u(\vartheta)\,d\vartheta,
\qquad
\beta\in C\comp^\infty([\frac{1}{2},2]).
\tag{\sec1.2}
$$

We assume that
the corank of the mixed Hessian $S\sb{x\thvar}$ 
in $(\sec1.1)$
(and hence the dimension of kernels of $d\pi\ssb{L}$,
$d\pi\ssb{R}$) is at most 1.

\proclaim{Terminology}
If the map $\pi\ssb{L}$ is of type at most $l$,
then we will say 
that the operator $T\sb{\lambda}$ has 
a singularity of type at most $l$ on the left.

Similarly with the singularity on the right.
\endproclaim

\proclaim{Theorem \sec1.1}
Let $T\sb{\lambda}$ be an oscillatory integral operator
of the form $(\sec1.2)$.
We assume that
$\rank S\sb{x\vartheta}\ge n-1$
and that $T\sb{\lambda}$ has singularities
of type at most $l$ on the left 
and at most $r$ on the right;
we denote $k=\min(l,r)$, $K=\max(l,r)$.

There are the following estimates
on the $L^2\to L^2$ action of $T\sb{\lambda}\sp{\hvar}$:
$$
\align
&\norm{T\sb{\lambda}\sp{\hvar}}
\le \const\lambda^{-\frac{n}{2}}
\hvar^{-\frac{1}{2}},
\qquad
\hvar\ge\lambda^{-\frac{1}{3}},
\tag{\sec1.3}
\\
&\norm{T\sb{\lambda}\sp{\hvar}}
\le \const\lambda^{-\frac{n+1}{2}}
\hvar^{-2},
\qquad
\lambda^{-(2+\frac{1}{k})^{-1}}
\le\hvar\le\lambda^{-\frac{1}{3}},
\tag{\sec1.4}
\\
&\norm{T\sb{\lambda}\sp{\hvar}}
\le \const\lambda^{-\frac{n}{2}}
\hvar^{-1+\frac{1}{2k}},
\qquad
\lambda^{-(2+\frac{1}{K})^{-1}}
\le\hvar\le\lambda^{-(2+\frac{1}{k})^{-1}},
\tag{\sec1.5}
\\
&\norm{T\sb{\lambda}\sp{\hvar}}
\le \const\lambda^{-\frac{n-1}{2}}
\hvar^{\frac{1}{2k}+\frac{1}{2K}},
\qquad
\lambda^{-\frac{1}{2}}
\le\hvar\le\lambda^{-(2+\frac{1}{K})^{-1}},
\tag{\sec1.6}
\endalign
$$
where the constants 
depend only 
on the bounds on derivatives of $\psi$ and $S$
in $(\sec1.2)$
(up to some finite order).

\endproclaim

Note that our methods are only applicable 
in the region $\hvar\ge\lambda^{-\frac{1}{2}}$
(in a certain sense, this is the restriction due to
the uncertainty principle).

\noindent
\remark{Remark}
The estimate $(\sec1.6)$ is also applicable 
to $T\sb{\lambda}$ localized to the region 
where 
$
\abs{h(x,\vartheta)}\le 2\hvar
$
(that is, when we no longer require 
$\abs{h(x,\vartheta)}\ge \hvar/2$).
We will denote such an operator by
$\bar{T}\sb{\lambda}\sp{\hvar}$:
$$
\bar{T}\sb{\lambda}\sp{\hvar}u(x)
=\underset{\R^n}\to\int
e^{i\lambda S(x,\vartheta)}
\psi(x,\vartheta)
\bar{\beta}(\hvar^{-1} h(x,\vartheta))
u(\vartheta)\,d\vartheta,
\qquad
\bar{\beta}\in C\comp^\infty([-2,2]).
\tag{\sec1.7}
$$
\endremark
Let the functions $\beta$ and
$\bar{\beta}$ satisfy
$
\sum\sb{\pm}\sum\sb{j=1}\sp{\infty}
\beta(\pm 2^{-j}t)+\bar{\beta}(t)=1,
$
for any $t\in\R$.
Then we can decompose $T\sb{\lambda}$
as
$$
T\sb{\lambda}=\sum\sb{\pm}\sum\sb{\hvar>\hvar\nut}
T\sb{\lambda}\sp{\pm\hvar}
+\bar{T}\sb{\lambda}\sp{\hvar\nut},
\qquad \hvar=2^{-N},\quad N\in\N,
\tag{\sec1.8}
$$
where the cut-off value
$\hvar\nut$ is to be chosen properly.
We can apply Theorem \sec1.1
to each $T\sb{\lambda}\sp{\hvar}$
(the estimates $(\sec1.3)$-$(\sec1.5)$)
and to $\bar{T}\sb{\lambda}\sp{\hvar}$ 
(the estimate $(\sec1.6)$).
This gives the following
estimate on $T\sb{\lambda}$:

\proclaim{Corollary 1}
Under the assumptions of Theorem \sec1.1,
$$
\norm{T\sb{\lambda}}
\le 
\const\lambda^{-\frac{n}{2}+\sup\sb{p}\delta},
\tag{\sec1.9}
$$
where
``the loss in the rate of decay at a point $p$''
is given by
$$
\delta(l,r)
=\frac{1}{2}\left(1-\frac{1}{2\min(l,r)}\right)
\left(1+\frac{1}{2\max(l,r)}\right)^{-1},
$$
with
$l$ and $r$ being the types of 
$\pi\ssb{L}$ and $\pi\ssb{R}$ at $p$.
The supremum in $(\sec1.9)$ is taken over
all points of ${\Cal C}$.
There is certainly no loss of smoothness
at non-singular points:
we define $\delta(0,0)=0$.

\endproclaim

This is weaker (except when $l=1$ or $r=1$)
than the optimal result
(proved in \cite{PhSt\ye{97}} for $n=1$)
which we might expect:
$
\delta\sb{\text{opt}}(l,r)
=\frac{1}{2}\left(1+\frac{1}{l}+\frac{1}{r}\right)^{-1}.
$

\subhead
Damping estimates
\endsubhead
We can use Theorem \sec1.1 for deriving
the damping estimates.
According to the estimates
$(\sec1.3)$--$(\sec1.6)$,
the series
$
\sum\sb{\pm}\sum\sb{\hvar>\lambda^{-\frac{1}{2}}}
\hvar\norm{T\sb{\lambda}\sp{\pm\hvar}}
+\lambda^{-\frac{1}{2}}\norm{
 \bar{T}\sb{\lambda}\sp{\lambda^{-\frac{1}{2}}}
}
$
(where $\hvar$ varies dyadically, 
as in $(\sec1.8)$),
is bounded by $\const\lambda^{-\frac{n}{2}}$.
Hence, if $U\sb{\lambda}$ is an operator like
$(\sec1.1)$ but with a damping factor 
of magnitude $\le\const\abs{\det S\sb{x\thvar}}$,
then
$
\norm{U\sb{\lambda}}
\le\const\lambda^{-n/2}.
$
This proves the following result:

\proclaim{Corollary 2}
Let $U\sb{\lambda}$ be a compactly supported 
oscillatory integral operator
of the form $(\sec1.1)$
with singularities of finite type
on both sides.
If the density $\psi$ vanishes on the critical variety 
$\varSigma
=\{(x,\thvar)\,\,\vert\,\, \det S\sb{x\thvar}(x,\thvar)=0\}$
so that
$
\abs{\psi}\le\const\abs{\det S\sb{x\thvar}},
$
then $U\sb{\lambda}$
has the same decay of its $L^p\rightarrow L^p$
norm as non-singular oscillatory integral operators:
$$
\norm{U\sb{\lambda}}\sb{L^p\rightarrow L^p}
\le \const\lambda^{-\frac{n}{2}+n\abs{\frac{1}{p}-\frac{1}{2}}},
\qquad 1\le p\le\infty.
\tag{\sec1.10}
$$

\endproclaim

Note that we have interpolated the $L^2$ estimates
with the trivial $L^1$ and $L^\infty$ estimates
(which are uniform in $\lambda$).

According to \cite{GrSe\ye{94}},
the $L^2$ estimate in
Corollary 2
implies the analogous result for Fourier integral operators:

\proclaim{Corollary 3}
Let $A\in I^m(X, Y,{\Cal C})$
be a Fourier integral operator
associated to 
a canonical relation ${\Cal C}$
such that
the projections $\pi\ssb{L}$, $\pi\ssb{R}$
are of corank at most 1
and have finite types everywhere.
If the symbol of $A$ vanishes on the critical 
variety of  $\pi\ssb{L}$ and $\pi\ssb{R}$, 
$
\abs{\sigma(A)}\le\const\abs{\det d\pi},
$
$d\pi$ being the Jacobi matrix
of either of $\pi\ssb{L}$, $\pi\ssb{R}$,
then for any real $s$
$$
A:\;
H\comp^s(Y)\rightarrow H\loc^{s-m}(X).
$$
\endproclaim

\subhead
$L^p$ estimates
\endsubhead
Let us say a few words about 
$L^p\rightarrow L^p$ estimates
on oscillatory integral operators.
They can be derived
by interpolating $L^2$ estimates with
$L^1\rightarrow L^1$ and $L^\infty\rightarrow L^\infty$ 
estimates:

\proclaim{Theorem \sec1.2}
Let $\rank S\sb{x\thvar}\ge n-1$.
If $T\sb{\lambda}$ 
has singularities of type at most
$l$ on the left and at most $r$ on the right,
then
$$
\align
&\norm{T\sb{\lambda}\sp{\pm\hvar}}\sb{L^1\rightarrow L^1}
\le\const
\hvar^{\frac{1}{r}},
\tag{\sec1.11}
\\ \\
&\norm{T\sb{\lambda}\sp{\pm\hvar}}\sb{L^\infty\rightarrow L^\infty}
\le\const
\hvar^{\frac{1}{l}}.
\tag{\sec1.12}
\endalign
$$

The same estimates 
are satisfied for $\bar{T}\sb{\lambda}\sp{\hvar}$.
\endproclaim

We may apply Theorems \sec1.1, \sec1.2
to derive the 
$L^p\rightarrow L^p$ estimates on 
$T\sb{\lambda}$.
Both $L^2$ and $L^1$, $L^\infty$ estimates 
on $\bar{T}\sb{\lambda}\sp{\hvar}$ 
become better for smaller values of $\hvar$;
by interpolation, we see that this is also
true for $L^p$ estimates for any $1\le p\le \infty$.
The estimates 
on $T\sb{\lambda}\sp{\hvar}$ have
a more complicated behavior:
$L^1$, $L^\infty$ estimates
become better for smaller values of $\hvar$,
while 
$L^2$ estimates ``blow up'' as $\hvar\rightarrow 0$.
Therefore, $L^p\rightarrow L^p$ estimates
on $T\sb{\lambda}\sp{\hvar}$
improve as $\hvar\rightarrow 0$
only if $p$ is outside a certain neighborhood of $p=2$.
In this case, the norms on the operators
in the dyadic decomposition $(\sec1.8)$
only become better as $\hvar$ becomes smaller,
and 
we conclude that 
the estimate on the entire $T\sb{\lambda}$ 
is determined by operators which are truncated off
the critical variety (large values of $\hvar$)
and hence coincides with the norm of non-degenerate
oscillatory integral operators.

In a certain neighborhood of $p=2$,
we need to glue the diverging $L^p$ estimates 
on $T\sb{\lambda}\sp{\hvar}$ 
with the estimate on $\bar{T}\sb{\lambda}\sp{\hvar\nut}$,
at some point $\hvar\nut>\lambda^{-\frac{1}{2}}$.

At some ``boundary values'' of $p$,
the $L^p$ estimates on $T\sb{\lambda}\sp{\hvar}$
are neither improving nor blowing up when $\hvar$
becomes small. 
Therefore, all the terms in $(\sec1.8)$
have the same bounds, 
and we are getting a factor 
$\ln \hvar\nut^{-1}\sim\ln\lambda$
(this is the number of terms in $(\sec1.8)$).

The $L^p$ estimates we obtain in this fashion 
are optimal only if the canonical relation
associated to $T\sb{\lambda}$
has a Whitney fold at least on one side
(this is when we know the optimal $L^2\to L^2$ estimates
\cite{Co\ye{98}}):

\proclaim{Corollary 4}
Let $T\sb{\lambda}$
be a compactly supported 
oscillatory integral operator 
with a fold singularity on the left.
If the singularity on the right is of type at most 
$r$,
then
the operator $T\sb{\lambda}$ has the same continuity
properties in $L^p$,
for $p<\frac{r+2}{r+1}$ and for $p>3$,
as a non-singular oscillatory integral operator:
$$
\norm{T\sb{\lambda}}\sb{L^p\rightarrow L^p}
\le\const
\lambda^{-\frac{n}{2}+n\abs{\frac{1}{p}-\frac{1}{2}}},
\qquad 1\le p <\frac{r+2}{r+1},
\quad 3< p\le\infty.
\tag{\sec1.13}
$$
For $\frac{r+2}{r+1}\le p\le 3$
the estimates are obtained by the interpolation with the 
$L^2$ estimates,
$$
\norm{T\sb{\lambda}}\sb{L^p\rightarrow L^p}
\le\const
\lambda^{-\frac{n}{2}+\left(4+\frac{2}{r}\right)^{-1}}.
\tag{\sec1.14}
$$
These estimates are sharp for $\frac{r+2}{r+1}<p<3$.
At the endpoints $p=\frac{r+2}{r+1}$ and $p=3$,
we can only prove weak estimates
(with the extra factor $\ln \lambda$).
\endproclaim

\demo{Remark}
The estimates in Theorem \sec1.2
may be improved if certain additional conditions
on the projections are satisfied, and this
in turn
leads to the estimate $(\sec1.13)$ in Corollary 4
to be true for a wider range of values of $p$.
For example, 
if $d\sb{x}(\det S\sb{x\thvar})\at{\varSigma}\ne 0$
(this condition is satisfied
if $\pi\ssb{R}$ has 
a strong $\sing{1\sb{r}}\sp{+}$-singularity, 
in the sense of \cite{GrSe\ye{97a}}),
then 
$\norm{T\sb{\lambda}\sp{\hvar}}\sb{L^1 \to L^1}
\le\const\hvar$,
and then one can easily prove that the estimate
$(\sec1.13)$ 
is valid in the range 
$1\le p<\frac{3}{2}$ and $3<p\le\infty$.
We will not discuss this issue here.
\enddemo

\subhead
$L^p$-$L^q$ estimates
\endsubhead
There are two more estimates
which hold for both 
$T\sb{\lambda}\sp{\pm\hvar}$ 
and $\bar{T}\sb{\lambda}\sp{\hvar}$:
$$
\norm{T\sb{\lambda}\sp{\pm\hvar}}
\sb{L^1\rightarrow L^\infty}
\le\const,
$$
which is trivially satisfied,
and 
$$
\norm{T\sb{\lambda}\sp{\pm\hvar}}
\sb{L^\infty\rightarrow L^1}
\le \const\hvar,
$$
which is satisfied if
$d\sb{x,\thvar} (\det S\sb{x\thvar})\at{\varSigma}\ne 0$.
The interpolation 
yields a variety of $L^p$-$L^q$ estimates
on $T\sb{\lambda}$
which we do not discuss.
For the case of oscillatory integral operators
with two-sided Whitney folds, see \cite{GrSe\ye{97b}}.

We will prove Theorem \sec1.2 in Section \sec2
and Theorem \sec1.1 in Sections \sec3, \sec4, and \sec5.

\head 
\sec2. Using finite type conditions:
$L^1$ and $L^\infty$ estimates
\endhead

We illustrate 
how the finite type conditions work
on the simplest example:
we derive
$L^1$ and $L^\infty$ estimates
formulated in Theorem \sec1.2.
We need certain preparation:
We choose local coordinates 
$x=(x\sp{\prime},x_n)$ and $\thvar=(\thvar\sp{\prime},\thvar_n)$
so that $S\sb{x\cpr \thvar\cpr}$ is non-degenerate
(recall that the corank of $S\sb{x\thvar}$ is at most 1),
and consider the map
$\pi\ssb{R}\at{\vartheta}$ as a composition
$$
(x\sp{\prime},x_n)
\overset{
  \pi\sp{\prime}
}
\to\longmapsto
(\eta\sp{\prime}=S\sb{\vartheta\cpr}(x,\vartheta),x_n)
\overset{
 \pi\sp{s}
}\to\longmapsto
(\eta\sp{\prime},\eta_n=S\sb{\vartheta\csb{n}}).
\tag{\sec2.1}
$$
According to the condition of the theorem 
that $\pi\ssb{R}$ is of type at most $r$,
we may assume that
$$
\vect{K}\ssb{R}^{r}
h\ge\varkappa\ssb{R}>0,
\tag{\sec2.2}
$$
where
the vector field
$\vect{K}\ssb{R}
=\left(\p\sb{x\csb{n}}
\right)\sb{\eta\cpr=S\csb{\thvar\cpr}\text{ fixed}}
$
has the property
$\vect{K}\ssb{R}\at{h(x,\thvar)=0}\in\ker d\pi\ssb{R}$.
Its explicit form is
$$
\vect{K}\ssb{R}
=\left(\p\sb{x\csb{n}}
\right)\sb{\eta\cpr=S\csb{\thvar\cpr}\text{ fixed}}
=\p\sb{x\csb{n}}
-S\sp{\vartheta\cpr x\cpr}S\sb{x\sb{n}\vartheta\cpr}
\p\sb{x\cpr},
\tag{\sec2.3}
$$
where
$S\sp{\vartheta\cpr x\cpr}(x,\thvar)$
is the inverse to the matrix $S\sb{x\cpr\vartheta\cpr}$ 
at a point $(x,\thvar)$.

\demo{Proof of Theorem \sec1.2}
We will write the generic notation $\beta$ for the
localizing functions $\beta$, $\bar{\beta}$;
the argument is the same
for both 
$T\sb{\lambda}\sp{\pm\hvar}$ 
and $\bar{T}\sb{\lambda}\sp{\hvar}$.
The key property of these operators 
is the small size of the support
of their integral kernels
``in the critical direction'';
we are to estimate this size using the finite type
conditions.

We have:
$$
\align
\norm{T\sb{\lambda}\sp{\hvar}u}\sb{L^1}
&\le\iint dx\,d\thvar\,
\abs{\psi(x,\thvar)\beta(\hvar^{-1}h(x,\thvar))
u(\thvar)}
\\
&\le
\norm{u}\sb{L^1}
\cdot\const
\underset{\thvar}\to\sup
\int dx\, 
\abs{\psi(x,\thvar)
\beta(\hvar^{-1}h(x,\thvar))}.
\endalign
$$
The bound which we need for the proof of $(\sec1.11)$,
$$
\int dx
\abs{\psi(x,\thvar)
 \beta(\hvar^{-1}h(x,\thvar))}
\le \const\hvar^{\frac{1}{r}},
\tag{\sec2.4}
$$
is due to the assumptions that 
the map $(x,\thvar)\mapsto(\thvar,S\sb{\thvar})$
is of type $r$.
We change the variables of integration to
$\eta\sp{\prime}=S\sb{\thvar\cpr}(x,\thvar)$
and $t=x_n$:
$$
\int dx\abs{\psi(x,\thvar)\beta(\hvar^{-1}h(x,\thvar))}
=
\int
\frac{d\eta\sp{\prime}}{\abs{\det S\sb{x\cpr \thvar\cpr}}}\,dt
\abs{\psi\beta(\hvar^{-1}h)}.
$$
We claim that the integration with respect to $t$ 
contributes $\const\hvar^{\frac{1}{r}}$
(while the integration with respect to $\eta\sp{\prime}$
is over the compact domain).
Indeed, 
since 
$\p\sb{t}=\left(\p\sb{x\csb{n}}\right)\sb{\eta\cpr}
=\vect{K}\ssb{R}$,
we know from $(\sec2.2)$
that
$
\p\sb{t}^{r}h\ge\varkappa\ssb{R}>0.
$
Now everything follows from the 
following lemma:

\proclaim{Lemma \sec2.1}
Let $h\in C^{r}(\R)$ be a function 
such that 
$\abs{h^{(r)}(t)}\ge\varkappa>0$ 
for $t$ in some interval $I\subset\R$.
Then the set
$
I\sp{\hvar}=\{t\in I\,\,\vert\,\, \abs{h(t)}<\hvar\}
$
consists of at most $2^{r-1}$ intervals
$I\sp{\hvar}\sb{\sigma}$,
possibly with joint ends,
with each of them being of measure
$
\abs{I\sp{\hvar}\sb{\sigma}}
\le\left({2r!}/{\varkappa}\right)^{\frac{1}{r}}
\hvar^{\frac{1}{r}}.
$

\endproclaim

This lemma is well-known; see, e.g., \cite{Ch\ye{85}}.
Let us give a proof which also motivates
the partition of 1 which will follow in Section \sec3.

\demo{Proof}
First, we take $\sigma$ to be a set
of $r-1$ signs, $\sigma=(\sigma_1,\,\dots,\,\sigma_{r-1})$,
$\sigma_j=\pm 1$.
We define
$$
I\sb{\sigma}=\{t\in I\,\,\vert\,\,
\sigma_j h^{(j)}(t)\ge 0,\,j=1,\,\dots,\,r-1 \}.
$$
Clearly, 
$I=\cup\sb{\sigma} I\sb{\sigma}$.
Since $h^{(r)}$ does not change sign, 
$h^{(r-1)}$ is monotone and hence the set
$\{t\in I\,\,\vert\,\, \sigma\sb{r-1}h^{(r-1)}\ge 0\}$
is connected.
Continuing by induction, we conclude that
$I\sb{\sigma}$ is also connected.
We now define 
$
I\sp{\hvar}\sb{\sigma}=I\sp{\hvar}\cap I\sb{\sigma},
$
which is also connected
(since $h$ is monotone on each $I\sb{\sigma}$).
We parameterize 
$I\sp{\hvar}\sb{\sigma}$ by $t$, 
so that
$t$ changes from 
$0$ to $\delta\equiv\abs{I\sb{\sigma}\sp{\hvar}}$.
Then for $0\le t\le \delta$ we have:
$$
\text{either}
\quad
\sigma_{r-1}h^{(r-1)}(t)\ge \varkappa t
\quad\text{or}\quad
\sigma_{r-1}h^{(r-1)}(t)\ge \varkappa {\delta-t}.
$$
The rest is by induction;
we will arrive at
$$
\text{either}
\quad
\sigma_{1}h^{\prime}(t)
\ge \frac{\varkappa}{(r-1)!}t^{r-1}
\quad\text{or}\quad
\sigma_{1}h^{\prime}(t)
\ge \frac{\varkappa}{(r-1)!}(\delta-t)^{r-1},
$$
and hence
$
\abs{h(\delta)-h(0)}
\ge\frac{\varkappa}{r!}\delta^{r}.
$
The {\it a priori} bound $\abs{h(\delta)-h(0)}<2\hvar$ 
gives the desired estimate on 
$\delta=\abs{I\sb{\sigma}\sp{\hvar}}$.
\qed
\enddemo

For the $L^\infty\to L^\infty$ estimate,
we derive
$$
\align
\norm{T\sb{\lambda}\sp{\hvar}u}\sb{L^\infty}
&\le\underset{x}\to\sup
\int d\thvar\,
\abs{\psi(x,\thvar)
\beta(\hvar^{-1}h(x,\thvar))
u(\thvar)}
\\
&\le
\norm{u}\sb{L^\infty}
\cdot\const
\underset{x}\to\sup
\int d\thvar\, 
\abs{\psi(x,\thvar)
\beta(\hvar^{-1}h(x,\thvar))}.
\endalign
$$
As above, we may prove that
if the map $(x,\thvar)\mapsto(x,S\sb{x})$
is of type $l$, then
$\int d\thvar\abs{\psi\beta(\hvar^{-1}h)} 
\le \const\hvar^{\frac{1}{l}}$,
and the estimate $(\sec1.12)$ follows.
This completes the proof of Theorem \sec1.2
\qed
\enddemo

\head \sec3. Localizations
\endhead

We are now going to prove Theorem \sec1.1.
Thus, let both $\pi\ssb{L}$
and $\pi\ssb{R}$ be of corank at most $1$
and have finite types.
We assume that 
on the support of the integral kernel of $T\sb{\lambda}$
the types of the projections
are at most $l$ and $r$, respectively.
Our statements for the cases $l=1$ or $r=1$
already follow from \cite{Co\ye{97}},
so we assume that $l>1$ and $r>1$.
For the definiteness, we will also assume that 
$l\ge r$.

The argument is the same for 
$T\sb{\lambda}\sp{\pm\hvar}$;
for our convenience we will always consider
$T\sb{\lambda}\sp{\hvar}$
(that is, we always assume that
$\det S\sb{x\thvar}$ is positive).

We split
the integral kernel of $T\sb{\lambda}\sp{\hvar}$
into pieces,
in the spirit of the proof of Lemma \sec2.1.
For this, we 
pick a smooth function $\rho$,
$\supp\rho\subset[-1,\infty)$,
such that $\rho(t)+\rho(-t)=1$,
and introduce 
the following partition of $1$:
$$
\aligned
&1=\underset{\sigma}\to\sum 
\rho\sb{\sigma}\sp{\hvar}(x,\vartheta),
\qquad 
\sigma=(\sigma\sb{1},\,\dots,\,\sigma\sb{l-1}),
\quad
\sigma\sb{j}=\pm 1,
\\
&\rho\sp{\hvar}\sb{\sigma}(x,\vartheta)
=\prod\sb{j=1}^{l-1}
\rho(\hvar^{-1}\sigma\sb{j}
\vect{K}\ssb{R}^{j}h(x,\vartheta)).
\endaligned
\tag{\sec3.1}
$$
Here the vector field
$
\vect{K}\ssb{R}
=\p\sb{x\csb{n}}
-S\sp{\vartheta\cpr x\cpr}S\sb{x\sb{n}\vartheta\cpr}
\p\sb{x\cpr}
$
is the same as in $(\sec2.3)$.

Analogously, we introduce the partition
$$
\aligned
&1=\underset{\varsigma}\to\sum 
\varrho\sb{\varsigma}\sp{\hvar}(x,\vartheta),
\qquad 
\varsigma=(\varsigma\sb{1},\,\dots,\,\varsigma\sb{r-1}),
\quad
\varsigma\sb{j}=\pm 1,
\\
&\varrho\sp{\hvar}\sb{\varsigma}(x,\vartheta)
=\prod\sb{j=1}^{r-1}
\rho(\hvar^{-1}\varsigma\sb{j}
\vect{K}\ssb{L}^{j}h(x,\vartheta)),
\qquad
\vect{K}\ssb{L}
=\p\sb{\vartheta\csb{n}}
-S\sp{\vartheta\cpr x\cpr}S\sb{x\cpr\vartheta\csb{n}}
\p\sb{\vartheta\cpr}.
\endaligned
\tag{\sec3.2}
$$
Of course,
$
\vect{K}\ssb{L}\at{h(x,\vartheta)=0}
\in\ker d\pi\ssb{L}.
$

Note that
$(i)$
the summation indexes $\sigma$ and $\varsigma$
in $(\sec3.1)$, $(\sec3.2)$
take finitely many values,
and 
$(ii)$ these are {\it admissible} partitions,
in the sense that
$$
\abs{\p\sb{x}\sp{\alpha}\p\sb{\vartheta}\sp{\beta}
\rho\sb{\sigma}\sp{\hvar}(x,\vartheta)}
\le C\sb{\alpha\beta}\hvar^{-\abs{\alpha}-\abs{\beta}},
$$
so that only $\hvar^{-1}$ can be contributed
during integrations by parts 
which will follow later in the argument.
We continue the proof individually for each
of the pieces of $T\sb{\lambda}\sp{\hvar}$
with fixed $\sigma$, $\varsigma$.

We use the ``fine'' partitions of $1$,
$$
\aligned
&1=\underset{X\in\Z^n}\to\sum
\chi(\hvar^{-1}\eta\sp{\prime}-X\sp{\prime})
\chi(\hvar^{-1}x_n-X_n),
\qquad \eta\sp{\prime}\equiv S\sb{\vartheta\cpr}(x,\vartheta),
\\
&1=\underset{\varTheta\in\Z^n}\to\sum
\chi(\hvar^{-1}\xi\sp{\prime}-\varTheta\sp{\prime})
\chi(\hvar^{-1}\vartheta_n-\varTheta_n),
\qquad \xi\sp{\prime}\equiv S\sb{x\cpr}(x,\vartheta),
\endaligned
\tag{\sec3.3}
$$
where $\chi$ is a certain smooth function 
supported in the unit ball in $\R^n$.
Multiplying the integral kernel of the operator
$T\sb{\lambda}\sp{\hvar}$
by the above functions,
we decompose $T\sb{\lambda}\sp{\hvar}$ into
$
T\sb{\lambda}\sp{\hvar}=
\underset{X,\varTheta\in\Z^n}\to\sum
\left(T\sb{\lambda}\sp{\hvar}\right)\sb{X\varTheta}.
$

\subhead{Convexity}
\endsubhead
We will use the fact that the map
$\pi\ssb{R}\at{\thvar}:\;x\mapsto S\sb{\thvar}(x,\thvar)$ 
(and similarly $\pi\ssb{L}\at{x}$)
satisfies certain convexity condition:
Given $\vartheta$, then for any $x$, $y$
on a connected set where ${\det S\sb{x\thvar}}\ge\hvar/2$
the following inequality holds:
$$
\abs{S\sb{\vartheta}(x,\vartheta)-S\sb{\vartheta}(y,\vartheta)}
\ge\const\hvar\abs{x-y}.
\tag{\sec3.4}
$$
Let us show
(sketching the argument from \cite{Co\ye{97}})
that 
the property $(\sec3.4)$ 
is satisfied on the support of each of 
$\sigma,\,\varsigma$-pieces of $T\sb{\lambda}\sp{\hvar}$.
The map $\pi\sp{\prime}$
in $(\sec2.1)$ is a diffeomorphism and 
hence we may assume that
$
\abs{\eta\sp{\prime}(x)-\eta\sp{\prime}(y)}
\ge\const\abs{\det S\sb{x\cpr\thvar\cpr}}\cdot\abs{x-y}.
$
We now need to investigate the map
$\pi\sp{s}\at{\eta\sp{\prime}}:\;
x\sb{n}\mapsto \eta\sb{n}=S\sb{\thvar\csb{n}}(x,\vartheta)$.
Let us denote by ${\Cal L}$ 
the line segment from 
$(\eta\sp{\prime},x_n)$ to $(\eta\sp{\prime},y_n)$.
We have:
$$
\abs{\eta\sb{n}(\eta\sp{\prime},y_n)-\eta\sb{n}(\eta\sp{\prime},x_n)}
\ge\abs{y\sb{n}-x\sb{n}}
\cdot\underset{{\Cal L}}\to\inf
\abs{\left(\p\sb{x\csb{n}}\right)\sb{\eta\sp{\prime}}
\!\eta\sb{n}}
\tag{\sec3.5}
$$
We need to show that the factor at 
$\abs{y\sb{n}-x\sb{n}}$
in the right-hand side of $(\sec3.5)$
is of magnitude $\hvar$,
and then
the inequality $(\sec3.4)$ follows.

The value of the derivative
$\left(\p\sb{x\csb{n}}\right)\sb{\eta\sp{\prime}}\!\eta\sb{n}$
can be determined from 
the decomposition
$\pi\ssb{R}\at{\vartheta}
=\pi\ssb{R}\sp{s}\at{\thvar}
\circ\pi\ssb{R}\sp{\prime}\at{\thvar}$.
Considering the determinants 
of the Jacobi matrices in $(\sec2.1)$,
$
J(\pi\ssb{R}\sp{s}\at{\thvar})
\cdot J(\pi\ssb{R}\sp{\prime}\at{\thvar})
=J(\pi\ssb{R}\at{\vartheta}),
$
we obtain
$
\left(\p\sb{x\csb{n}}\right)\csb{\eta\cpr}\!\eta\sb{n}
\cdot\det S\sb{x\cpr\vartheta\cpr}
=h(x,\vartheta).
$
Hence,

\proclaim{Lemma \sec3.1}
There is the relation
$
\left(\p\sb{x\csb{n}}\right)\sb{\eta\sp{\prime}}
\!\eta\sb{n}
=\frac{h(x,\vartheta)}{\det S\sb{x\cpr\vartheta\cpr}}.
$
\endproclaim

We also need to check that
$\ h(x,\vartheta)\ge\const\hvar$
on a line between $x_n$ and $y_n$:

\proclaim{Lemma \sec3.2}
If $\abs{{\Cal L}}\le\const$,
then $h\ge\frac{\hvar}{4}$
everywhere on ${\Cal L}$.

\endproclaim

Hence, we admit that the line segment
${\Cal L}$
could be not entirely on the support 
of the integral kernel of $T\sp{\hvar}$,
where $h\ge\hvar/2$.

\demo{Proof}
Let $t$ be a parameter on the line segment ${\Cal L}$,
changing from $t=0$ at the point $(\eta\sp{\prime},x_n)$
to $t=\abs{x_n-y_n}$ at the point 
$(\eta\sp{\prime},y_n)$.
We consider $h(x,\vartheta)\at{{\Cal L}}$ 
as a function of $t$;
$\p_t=\left(\p\sb{x\csb{n}}\right)\sb{\eta\sp{\prime}}$.
Since both $(\eta\sp{\prime},x_n)$ and $(\eta\sp{\prime},y_n)$
are on the support of $\sigma,\,\varsigma$-piece
of $T\sb{\lambda}\sp{\hvar}$,
we know that 
at $t=0$ and at $t=\abs{x_n-y_n}$
$$
\align
&h\ge\hvar/2,
\tag{\sec3.6}
\\
&\sigma_{j}
\cdot (\p\sb{t})^{j}h 
\ge -\hvar,
\qquad \forall\,\,j\in\N,\quad j<r.
\tag{\sec3.7}
\endalign
$$
Due to the finite type conditions
on both projections from ${\Cal C}$,
we also know that
$
(\p\sb{t})^{r} h > 0
$
(or instead $<0$),
for all $t$ between $0$ and $\abs{x_n-y_n}$
(see $(\sec2.2)$).

If we assumed that in 
$({\sec3.6})$ and $({\sec3.7})$
$\hvar=0$, then we would conclude
by induction
that all $(\p\sb{t})^{j}h(t)$, $j<r$,
were monotone functions 
which did not change the signs
between $t=0$ and $t=\abs{x_n-y_n}$,
and hence $h(t)$ would be concluded monotone
(see the proof of Lemma \sec2.1).
Since $\hvar\ne 0$,
the above conclusion 
is true modulo the error of magnitude
$\hvar$; hence,
the function $h(t)$ is ``almost monotone'' (its derivative
is greater than $-\const\hvar$ or less than $\const\hvar$),
therefore the value of $h(t)$ can not drop
below $\hvar/4$
as long as $t$ is between $0$ and $\abs{x_n-y_n}$
and as long as at the boundary points
the value of $h$ is not less than $\hvar/2$.
We also need to assume that 
$\abs{x_n-y_n}$
is not too large.
\qed
\enddemo

\head
\sec4. Almost orthogonality relations for different pieces
\endhead

We are going to apply the Cotlar-Stein lemma on $L^2$
almost orthogonality \cite{St\ye{93}}.
For our convenience, let us formulate this result here.

\proclaim{Cotlar-Stein Lemma}
Let $E$ and $F$ be the Hilbert spaces,
and let
$
\{T_i\,\,\vert\,\, i\in\Z\}
$
be a family of continuous operators $E\to F$
which satisfy the following conditions:
$$
\norm{T\sb{i}^\ast T\sb{j}}\le  a(i,j),
\qquad\qquad
\norm{T\sb{i} T\sb{j}^\ast}\le  b(i,j),
\tag\sec4.1
$$
where $a(i,j)$ and $b(i,j)$ are non-negative
functions on $\Z\times\Z$.
If $a$ and $b$ satisfy
$$
A
\equiv
\underset{i}\to\sup\sum\sb{j}  a^{\frac{1}{2}}(i,j)
<\infty,
\qquad
B\equiv
\underset{i}\to\sup\sum\sb{j}  b^{\frac{1}{2}}(i,j)
<\infty,
\tag\sec4.2
$$
then the formal sum $\sum_i T_i$
converges (in the weak operator topology)
to a continuous operator 
$
\;T:\,E\rightarrow F,\;
$ 
which is bounded by
$$
\norm{T}
\le A^{\frac{1}{2}}B^{\frac{1}{2}}.
\tag\sec4.3
$$

\endproclaim

The details of the proof are in \cite{St\ye{93}}.

Now we are going to investigate the almost orthogonality
relations for the operators
$\left(T\sb\lambda\sp\hvar\right)\sb{X \varTheta}$,
with respect to different multi-indexes $\{X \varTheta\}$.

\subhead{Almost orthogonality
with respect to different $\varTheta$, $W$}
\endsubhead
Let us consider 
the behavior of the compositions
$\left(T\sb{\lambda}\sp{\hvar}\right)\sb{X\varTheta}
\left(T\sb{\lambda}\sp{\hvar}\right)\sb{Y W}^\ast$
with respect to different $\varTheta$ and $W$.
We will show that
if $X$ and $Y$ are fixed and if $W$
is also fixed,
then the composition 
$\left(T\sb{\lambda}\sp{\hvar}\right)\sb{X\varTheta}
\left(T\sb{\lambda}\sp{\hvar}\right)\sb{Y W}^\ast$
is different from zero only for finitely many values 
of $\varTheta$.
The integral kernel of such an operator is given by
$$
\aligned
K\sb{\varTheta W}(x,y)
&=\int_{\R^n}d\vartheta\,
\chi(\hvar^{-1}S\sb{x\cpr}(x,\vartheta)-\varTheta\sp{\prime})
\chi(\hvar^{-1}S\sb{y\cpr}(y,\vartheta)-W\sp{\prime})
\\
&\qquad\qquad
\chi(\hvar^{-1}S\sb{\vartheta\cpr}(x,\vartheta)-X\sp{\prime})
\chi(\hvar^{-1}S\sb{w\cpr}(y,\vartheta)-Y\sp{\prime})
\times\dots,
\endaligned
\tag{\sec4.4}
$$
and in addition we know that
$x_n\approx\hvar X_n$, $y_n\approx \hvar Y_n$, and
$\vartheta_n\approx\hvar\varTheta_n\approx\hvar W_n$.
(Each time, the error is at most $\hvar$.)
Let us consider the following system:
$$
\cases
&S\sb{y\cpr}(y\sp{\prime},y_n,\vartheta\sp{\prime},\vartheta_n)
=\hvar W\sp{\prime},
\\
&
S\sb{\vartheta\cpr}(y\sp{\prime},y_n,\vartheta\sp{\prime},\vartheta_n)
=\hvar Y\sp{\prime}.
\endcases
\tag{\sec4.5}
$$
Given $y_n$, $\thvar_n$, $Y\sp{\prime}$, and $W\sp{\prime}$,
we can solve this system for $y\sp{\prime}$ and $\vartheta\sp{\prime}$,
since the matrix
$\dsize
\frac{\p(\hvar W\sp{\prime},\hvar Y\sp{\prime})}
{\p(y\sp{\prime},\vartheta\sp{\prime})}
=\left[
\matrix 
S\sb{y\cpr\vartheta\cpr}
&S\sb{\vartheta\cpr\vartheta\cpr}
\\
S\sb{y\cpr y\cpr}
&S\sb{y\cpr\vartheta\cpr}
\endmatrix
\right]
$
is non-degenerate. 
(For this, we could have used 
certain preparation:
at some
point $(x\nut,\vartheta\nut)$ 
on the support of $T\sb{\lambda}$
we make
the matrices
$S\sb{x\cpr x\cpr}$ and $S\sb{\vartheta\cpr\vartheta\cpr}$
vanish, using the change
$S(x,\vartheta)\mapsto S(x,\vartheta)
+\Phi_1(x\sp{\prime})+\Phi_2(\vartheta\sp{\prime})$,
which is equivalent to 
a unitary transformation.
Possibly, we also use a restriction 
to a smaller neighborhood of $(x\nut,\thvar\nut)$.)
Now, since the parameters
$y_n$ and $\vartheta_n$
($\approx\hvar Y_n$ and $\approx\hvar\varTheta_n$)
and also the right-hand sides 
of the system $(\sec4.5)$
are determined 
with the error $\hvar$,
both $y\sp{\prime}$ and $\vartheta\sp{\prime}$
are determined with the error of the same magnitude.
Then, from $S\sb{\vartheta\sp{\prime}}(x\sp{\prime},x_n,\vartheta)
\approx\hvar X\sp{\prime}$,
$x_n\approx \hvar X_n$,
we determine $x\sp{\prime}$ (same error).
Hence, 
$\hvar \varTheta\sp{\prime}\approx S\sb{x\sp{\prime}}(x,\vartheta)$
is also determined with the error of magnitude $\hvar$.
We conclude that 
$\varTheta$ can take only finitely many values
(uniformly in $\lambda$, $\hvar$).
Note that the particular range of $\varTheta$ 
may depend on
the specific values of $X$, $Y$, and $W$.

\subhead{Almost orthogonality
with respect to different $X$, $Y$}
\endsubhead
The almost orthogonality of compositions
$
\left(T\sb{\lambda}\sp{\hvar}\right)\sb{X\varTheta}
\left(T\sb{\lambda}\sp{\hvar}\right)\sb{Y W}^\ast
$
with respect to different $X$ and $Y$
requires the integration by parts
in the expression for the integral kernel:
$$
\aligned
&K\left(
\left(T\sb{\lambda}\sp{\hvar}\right)\sb{X\varTheta}
\left(T\sb{\lambda}\sp{\hvar}\right)\sb{Y W}^\ast
\right)(x,y)
=\underset{\R^n}\to\int
e^{i\lambda(S(x,\vartheta)-S(y,\vartheta))}
\chi(\dots)
\beta(\dots)
\psi(\dots)
d^n\vartheta.
\endaligned
\tag{\sec4.6}
$$
Integration by parts
in the expression $(\sec4.6)$
shows that for any integer $N$
$$
\Abs{K\left(
\left(T\sb{\lambda}\sp{\hvar}\right)\sb{X\varTheta}
\left(T\sb{\lambda}\sp{\hvar}\right)\sb{Y W}^\ast
\right)(x,y)}
\le \const\sb{N}
\int 
\frac{d^n\thvar\,\chi(\dots)\beta(\dots)}
{\left[1+\lambda\hvar
\abs{S\sb{\thvar}(x,\thvar)-S\sb{\thvar}(y,\thvar)}
\right]^{2N}}.
\tag{\sec4.7}
$$
The factor $\hvar$ in the denominator
reflects the contribution of $\hvar^{-1}$
from each integration by parts
(to estimate the contribution of certain terms,
one needs to refer to $(\sec3.4)$).

We claim that
$$
\abs{S\sb{\thvar}(x,\thvar)-S\sb{\thvar}(y,\thvar)}
\ge
\const(\hvar\abs{X\sp{\prime}-Y\sp{\prime}}+\hvar^2\abs{X_n-Y_n}).
$$
For $\hvar\abs{X_n-Y_n}\le \abs{X\sp{\prime}-Y\sp{\prime}}$
this inequality is trivial, while for 
$\hvar\abs{X_n-Y_n}\ge \abs{X\sp{\prime}-Y\sp{\prime}}$
we use the 
convexity property $(\sec3.4)$ of $\pi\ssb{R}$:
$$
\abs{S\sb{\thvar}(x,\vartheta)-S\sb{\thvar}(y,\vartheta)}
\ge
\const\hvar\abs{x-y}
\ge\const\hvar^2\abs{X_n-Y_n}.
$$

Now we may rewrite the right-hand side of $(\sec4.7)$
as
$$
\frac{\const}
{\left[1+\lambda\hvar
(\hvar\abs{X\sp{\prime}-Y\sp{\prime}}+\hvar^2\abs{X_n-Y_n})
\right]^{N}}
\int 
\frac{d^n\thvar\,\chi(\dots)\beta(\dots)}
{\left[1+\lambda\hvar
\abs{S\sb{\thvar}(x,\thvar)-S\sb{\thvar}(y,\thvar)}
\right]^{N}}
$$
and apply the Schur lemma:
$$
\align
&
\Norm{
\left(T\sb{\lambda}\sp{\hvar}\right)\sb{X\varTheta}
\left(T\sb{\lambda}\sp{\hvar}\right)\sb{Y W}^\ast
}
\le
\int d^n x \Abs{K
\left(
\left(T\sb{\lambda}\sp{\hvar}\right)\sb{X\varTheta}
\left(T\sb{\lambda}\sp{\hvar}\right)\sb{Y W}^\ast
\right)
(x,y)}
\\
&\le
\frac{\const}
{\left[1+\lambda\hvar
(\hvar\abs{X\sp{\prime}-Y\sp{\prime}}+\hvar^2\abs{X_n-Y_n})
\right]^{N}}
\int 
\frac{d^n x \,d^n \thvar\,\chi(\dots)\beta(\dots)}
{\left[1+\lambda\hvar
\abs{S\sb{\thvar}(x,\thvar)-S\sb{\thvar}(y,\thvar)}
\right]^{N}}.
\endalign
$$

We will integrate in $x$ first.
If $\hvar\ge\lambda^{-\frac{1}{3}}$,
then
$$
\int 
\frac{d^n x}
{\left[1+\lambda\hvar
\abs{S\sb{\thvar}(x,\thvar)-S\sb{\thvar}(y,\thvar)}
\right]^{N}}
=\int 
\frac{d^n \{ S\sb{\thvar}\}}
{\abs{\det{S\sb{x\thvar}}}}
\cdot
\frac{1}
{\left[1+\lambda\hvar
\abs{S\sb{\thvar}(x,\thvar)-S\sb{\thvar}(y,\thvar)}
\right]^{N}}
$$
is bounded by
$\const(\lambda\hvar)^{-n}\hvar^{-1}.$
If instead $\hvar\le\lambda^{-\frac{1}{3}}$,
then a better bound is obtained 
when appealing the size of the support in $x_n$:
$$
\int\frac{d^n x \,\chi(\dots)\beta(\dots)}
{\left[1+\lambda\hvar
\abs{S\sb{\thvar}(x,\thvar)-S\sb{\thvar}(y,\thvar)}
\right]^{N}}
=\int
\frac{d^{n-1}\{S\sb{\thvar\cpr}\}}
{\abs{\det S\sb{x\cpr\thvar\cpr}}}
\cdot
\frac{dx_n\,\chi(\hvar^{-1}x_n-X_n)\dots }
{\left[1+\lambda\hvar
\abs{S\sb{\thvar}(x,\thvar)-S\sb{\thvar}(y,\thvar)}
\right]^{N}}.
$$
This expression is bounded by
$\le\const(\lambda\hvar)^{-n+1}\hvar$:
$(\lambda\hvar)^{-n+1}$ is due to the integration in 
$S\sb{\thvar\sp{\prime}}$, 
and $\hvar$ is due to the integration in $x_n$.

The integration with respect to $\thvar$
is performed as follows:
$$
\int d\thvar^n
\chi(\dots)
=\int\frac{d^{n-1}\{S\sb{y\cpr}\}\,d\thvar_n}
{\abs{\det S\sb{y\cpr\thvar\cpr}}}
\chi(\hvar^{-1}S\sb{y\cpr}(y,\thvar)-\varTheta\sp{\prime})
\chi(\hvar^{-1}\thvar_n-\varTheta_n)
\le\const\hvar^{n}.
$$
We conclude that
$$
\norm{
\left(T\sb{\lambda}\sp{\hvar}\right)\sb{X\varTheta}
\left(T\sb{\lambda}\sp{\hvar}\right)\sb{Y W}\sp{\ast}}
\le \const \tau^2
\left[1+\lambda\hvar^2\left(
\abs{X\sp{\prime}-Y\sp{\prime}}+
\hvar\abs{X_n-Y_n}
\right)\right]^{-N},
\tag{\sec4.8}
$$
where 
$\tau^2=\min(\lambda^{-n}\hvar^{-1},\lambda^{-n+1}\hvar^2)$.
One may think of $\tau$ 
as of the $L^2$-estimate on a generic operator 
$\left(T\sb{\lambda}\sp{\hvar}\right)\sb{X\varTheta}$.

Let us rewrite $(\sec4.8)$ as
$$
\norm{\left(T\sb{\lambda}\sp{\hvar}\right)\sb{X\varTheta}
\left(T\sb{\lambda}\sp{\hvar}\right)\sb{Y W}\sp{\ast}}
\le
\tau^2
a(X,\varTheta;Y,W),
\tag{\sec4.9}
$$
where $a(X,\varTheta;Y,W)\in C(\Z^{4n};\R\sb{+})$ is some 
function which ``measures'' the orthogonality
of operators.
Similarly to $(\sec4.2)$,
we define
$$
A=\underset{Y,W}\to\sup\,
\sum\sb{X,\varTheta}a^\frac{1}{2}(X,\varTheta,Y,W);
\tag{\sec4.10}
$$
roughly, this is ``the number of the operators
which are not orthogonal''
(imagine that $a$ takes values $0$ and $1$ only).
We proved earlier that for fixed
$X$, $Y$, and $W$, the multi-index $\varTheta$
takes only finitely many values; 
therefore 
the summation with respect to $\varTheta$ 
in $(\sec4.10)$
is over a finite region in $\Z^n$ 
and only contributes some factor
which is uniform in $\hvar$ and $\lambda$.

Before we proceed to the analysis of the summation in $X$,
let us say a few words about the compositions 
of the form $T^\ast T$.
For such compositions, we have estimates
similar to $(\sec4.9)$:
$$
\norm{\left(T\sb{\lambda}\sp{\hvar}\right)\sb{X\varTheta}
\sp{\ast}
\left(T\sb{\lambda}\sp{\hvar}\right)\sb{Y W}}
\le C \tau^2
b(X,\varTheta;Y,W);
\tag{\sec4.11}
$$
we define
$$
B=\underset{Y,W}\to\sup\,
\sum\sb{X,\varTheta}b^\frac{1}{2}(X,\varTheta;Y,W).
\tag{\sec4.12}
$$
By the symmetry, we know that the summation 
with respect to $X$ in $(\sec4.12)$
is over a finite subset in $\Z^n$.

We analyze the summation in $X$ in four different cases:
$\hvar\ge\lambda^{-\frac{1}{3}}$, 
$\lambda^{-(2+\frac{1}{r})^{-1}}
\le\hvar\le\lambda^{-\frac{1}{3}}$,
$\lambda^{-(2+\frac{1}{l})^{-1}}
\le\hvar\le\lambda^{-(2+\frac{1}{r})^{-1}},
$
and 
$\lambda^{-\frac{1}{2}}\le\hvar
\le\lambda^{-(2+\frac{1}{r})^{-1}}$.

\subhead{The case $\hvar\ge\lambda^{-\frac{1}{3}}$}
\endsubhead
The simplest case is when $\hvar\ge\lambda^{-\frac{1}{3}}$;
then $(\sec4.8)$ decreases faster than 
$\tau^2$ times any power of $\Abs{X-Y}^{-1}$,
and the sum $\sum\sb{X}$ 
in $(\sec4.10)$
is bounded uniformly in $\hvar$, $\lambda$,
so that $A\le\const$.
By the symmetry, $B$ in $(\sec4.12)$ is also uniformly bounded.

According to the Cotlar-Stein lemma $(\sec4.3)$,
the bound on $T\sb{\lambda}\sp{\hvar}$
is given by 
$$
\tau=\min(\lambda^{-\frac{n-1}{2}}\hvar,
\lambda^{-\frac{n}{2}}\hvar^{-\frac{1}{2}}),
$$
which is 
the root of the common factor in
$(\sec4.9)$, $(\sec4.11)$,
times the geometric mean
of $A$ and $B$:
$$
\norm{T\sb{\lambda}\sp{\hvar}}
\le
\tau
\sqrt{AB}
\le\const\lambda^{-\frac{n}{2}}\hvar^{-\frac{1}{2}}\sqrt{AB}
\le\const\lambda^{-\frac{n}{2}}\hvar^{-\frac{1}{2}},
\qquad\hvar\ge\lambda^{-\frac{1}{3}}.
\tag{\sec4.13}
$$

\subhead{The case 
$\lambda^{-(2+\frac{1}{r})^{-1}}
\le\hvar\le\lambda^{-\frac{1}{3}}$}
\endsubhead
Again, for each $X$, the summation in $\varTheta$ in
$
A=\sum\sb{X,\varTheta}a^\frac{1}{2}(X,\varTheta;Y,W)
$
is over a finite set of multi-indices.
But, if $\lambda\hvar^3\le 1$, 
then 
it follows from $(\sec4.8)$ that
the summation 
with respect to $X$ in $(\sec4.10)$
contributes 
$$
\sum\sb{X}
[1
+\lambda\hvar^2\abs{X\sp{\prime}-Y\sp{\prime}}
+\lambda\hvar^3\abs{X_n-Y_n}]^{-N}
\le\const(\lambda\hvar^3)^{-1}.
$$
Note that 
the summation with respect to $X\sp{\prime}$ 
in $(\sec4.10)$
is fine
(contributes a factor uniform in $\lambda$, $\hvar$)
as long as $\lambda\hvar^2\ge 1$.
Our conclusion is that
``the number of non-orthogonal operators'' is
controlled by
$$
A=\sum\sb{X,\varTheta}a^\frac{1}{2}(X,\varTheta;Y,W)
\le\const(\lambda\hvar^3)^{-1}.
\tag{\sec4.14}
$$
Similarly,
$
B=\sum\sb{X,\varTheta}b^\frac{1}{2}(X,\varTheta;Y,W)
\le\const(\lambda\hvar^3)^{-1},
$
and $(\sec4.13)$ becomes
$$
\norm{T\sb{\lambda}\sp{\hvar}}\le
\const\lambda^{-\frac{n-1}{2}}\hvar
\sqrt{AB}
=\const\lambda^{-\frac{n+1}{2}}\hvar^{-2},
\qquad
\lambda^{-\frac{1}{2}}\le\hvar\le\lambda^{-\frac{1}{3}}.
\tag{\sec4.15}
$$

\subhead{The case 
$\lambda^{-(2+\frac{1}{l})^{-1}}
\le\hvar\le\lambda^{-(2+\frac{1}{r})^{-1}}$}
\endsubhead
The estimate $(\sec4.15)$ is clearly inadequate
for small values of $\hvar$:
its derivation is based on the assumption that
$\det S\sb{x\vartheta}\ne 0$, and as a consequence the
estimate blows up when $\hvar\rightarrow 0$
and does not allow to estimate the 
contribution of some tiny neighborhood 
of the critical variety $\{\det S\sb{x\vartheta}=0\}$.
Let us get another estimate on $A$, trying to
count ``non-orthogonal terms'' directly.
In $(\sec4.14)$, we evaluated the sum
assuming that the number of terms
with different $X_n$ is infinite,
while certainly $X_n$ takes at most $\const \hvar^{-1}$
values.
More than that, 
if the projection $\pi\ssb{R}$ is of type $r$,
then there are only $\const\hvar^{-1+\frac{1}{r}}$
terms with different $X\sb{n}$. 
This is because the size of
support of the integral kernel of $T\sb{\lambda}\sp{\hvar}$
is not only compact, but also bounded 
in certain critical direction
by $\hvar^{\frac{1}{r}}$
(this is gained by the methods which are
very much the same as in Section \sec2;
we will show this in more detail
in Section \sec5).
This leads to the following bound on $A$ in $(\sec4.10)$:
$$
A
\le\const\hvar^{-1+\frac{1}{r}}.
\tag\sec4.16
$$

Now we proceed to deriving the resulting bounds on 
$T\sb{\lambda}\sp{\hvar}$.
If $\hvar\le\lambda^{-(2+\frac{1}{r})^{-1}}$,
then $(\sec4.16)$ gives a better bound on $A$
than $(\sec4.14)$.
We assume that 
$\hvar\ge\lambda^{-(2+\frac{1}{l})^{-1}}$,
so that a proper bound on $B$
is still $(\lambda\hvar^3)^{-1}$.
Hence, 
we can rewrite $(\sec4.13)$ as
$$
\norm{T\sb{\lambda}\sp{\hvar}}
\le \const\lambda^{-\frac{n-1}{2}}\hvar
\sqrt{\hvar^{-1+\frac{1}{r}}(\lambda\hvar^3)^{-1}}
\le \const\lambda^{-\frac{n}{2}}\hvar^{-1+\frac{1}{2r}}.
\tag\sec4.17
$$

\subhead{The case 
$\lambda^{-\frac{1}{2}}\le\hvar\le\lambda^{-(2+\frac{1}{l})^{-1}}$}
\endsubhead
In this region,
the best bounds on both $A$ and $B$
are due to the finiteness of types of the projections
$\pi\ssb{L}$ and $\pi\ssb{R}$:
$
A\le\const\hvar^{-1+\frac{1}{r}},
$
$
B\le\const\hvar^{-1+\frac{1}{l}}.
$
This gives
$$
\norm{T\sb{\lambda}\sp{\hvar}}
\le \const\lambda^{-\frac{n-1}{2}}\hvar
\sqrt{\hvar^{-1+\frac{1}{r}}\hvar^{-1+\frac{1}{l}}}
\le \const\lambda^{-\frac{n-1}{2}}
\hvar^{\frac{1}{2l}+\frac{1}{2r}}.
\tag\sec4.18
$$
Since the derivation of 
bounds $A\le\const\hvar^{-1+\frac{1}{r}}$ 
and $B\le\const\hvar^{-1+\frac{1}{l}}$ 
does not appeal
to the inequality $\det S\sb{x\thvar}\ge\hvar/2$
(for details, see Section \sec5),
we conclude that $\bar{T}\sb{\lambda}\sp{\hvar}$
also satisfies the estimate $(\sec4.18)$.

\head
\sec5. Almost orthogonality and finite type conditions 
\endhead

We are left to prove the bound $(\sec4.16)$.
As in Section \sec2,
we consider the map
$\pi\ssb{R}\at{\vartheta}$ 
as decomposed into
$$
(x\sp{\prime},x_n)
\overset{\pi\sp{\prime}}\to\longmapsto
(\eta\sp{\prime}=S\sb{\vartheta\cpr}(x,\vartheta),x_n)
\overset{\pi\sp{s}}\to\longmapsto
(\eta\sp{\prime},\eta_n=S\sb{\vartheta\csb{n}}).
\tag{\sec5.1}
$$
Since $\pi\ssb{R}$ is of type at most $r$,
we may assume that on the support of 
the integral kernel of the operator $T\sb\lambda$
we have a uniform bound
$$
\left(\p\sb{x\csb{n}}\right)\sb{\eta\sp{\prime}}^{r}
h\ge\varkappa\ssb{R}>0.
\tag{\sec5.2}
$$
We work in the space $(\eta\sp{\prime},x_n)$.
Consider the line segment ${\Cal L}$ which 
connects the points
$(\eta\sp{\prime}(x,\vartheta),x_n)$
and
$(\eta\sp{\prime}(y,\vartheta),y_n)$.
There are two cases:

\noindent
$\bullet$
The line segment ${\Cal L}$ is {\it outside} the 
$c\hvar^{1-\frac{1}{r}}$-cone
of the direction
$\pm\left(\p\sb{x\csb{n}}\right)\sb{\eta\sp{\prime}}$
($c$ should be sufficiently small; see later);
this corresponds to 
$$
\abs{S\sb{\vartheta\cpr}(x,\vartheta)
-S\sb{\vartheta\cpr}(y,\vartheta)}
\ge c\hvar^{1-\frac{1}{r}}\abs{x_n-y_n}.
\tag{\sec5.3}
$$
We derive that
the gradient of the phase 
function,
$
\lambda(S\sb{\vartheta}(x,\vartheta)
-S\sb{\vartheta}(y,\vartheta)),
$
is bounded in the absolute value from below by
$$
\lambda
\abs{S\sb{\vartheta\cpr}(x,\vartheta)
-S\sb{\vartheta\cpr}(y,\vartheta)}
\ge
\frac{1}{2}
\lambda\left(
\abs{S\sb{\vartheta\cpr}(x,\vartheta)
-S\sb{\vartheta\cpr}(y,\vartheta)}
+
c\hvar^{1-\frac{1}{r}}\abs{x_n-y_n}\right).
$$

Similarly to
how we arrived at $(\sec4.8)$,
we derive that for any integer $N$
$$
\norm{
\left(T\sb{\lambda}\sp{\hvar}\right)\sb{X\varTheta}
\left(T\sb{\lambda}\sp{\hvar}\right)\sb{Y W}\sp{\ast}}
\le\const
\lambda^{-n+1}\hvar^{2}\left[1+\lambda\hvar^2\left(
\abs{X\sp{\prime}-Y\sp{\prime}}
+
c\hvar^{1-\frac{1}{r}}\abs{X_n-Y_n}
\right)\right]^{-N}.
$$
Since $\lambda\hvar^{2}\ge 1$,
the summation $\sum\sb{X_n}$ in $(\sec4.10)$
contributes at most 
$$
\sum\sb{X_n}\left[
1
+c\lambda\hvar^{3-\frac{1}{r}}\Abs{X_n-Y_n}
\right]^{-N}
\le\const c^{-1} \hvar^{-1+\frac{1}{r}},
$$
in an agreement with $(\sec4.16)$.
(Recall at this point 
that the summation with respect to
$\varTheta$ in $(\sec4.10)$
is over a bounded set in $\Z^n$
and that the summation with respect to $X\sp{\prime}$
converges, as long as $\hvar\ge \lambda^{-\frac{1}{2}}$.)

\noindent
$\bullet$
Now assume that the line segment ${\Cal L}$ 
from $(\eta\sp{\prime}(x),x_n)$ to $(\eta\sp{\prime}(y),y_n)$
is {\it inside} the 
$c\hvar^{1-\frac{1}{r}}$-cone
of the directions
$\pm\left(\p\sb{x\csb{n}}\right)\sb{\eta\sp{\prime}}$.
If $\abs{x-y}\ge\const\hvar^{\frac{1}{r}}$,
then, using $(\sec5.2)$,
we will show that
$$
\abs{h(x,\vartheta)-h(y,\vartheta)}
\sim
\varkappa\ssb{R}(\hvar\abs{{\Cal L}})\sp{r},
\tag{\sec5.4}
$$
where 
$
\abs{{\Cal L}}
\equiv\text{dist}\,[(\eta\sp{\prime}(x),x_n),\,(\eta\sp{\prime}(y),y_n)].
$
Since the left-hand side of $(\sec5.4)$ can not be greater than
$4\hvar$
(the value of $\abs{h}$
at both points 
$(x,\vartheta)$ and $(y,\vartheta)$
is bounded by $2\hvar$),
we gain
the bound $\abs{{\Cal L}}\le\const\hvar^{\frac{1}{r}}$.
This, together with 
$\abs{{\Cal L}}\ge\abs{x_n-y_n}\approx\hvar\Abs{X_n-Y_n}$,
yields the desired restriction
$$
\abs{X_n-Y_n}
\le\const\hvar^{-1+\frac{1}{r}},
\tag{\sec5.5}
$$
which again leads to $(\sec4.16)$.

The detailed proof of $(\sec5.5)$
is in \cite{Co\ye{98}};
for the reader's convenience, 
we give here the sketch.
Let $t$ be a parameter on the line segment
${\Cal L}$, 
which changes
from $t=0$ at $\pi\sp{\prime}(x)$ 
to $t=\abs{{\Cal L}}$ at 
$\pi\sp{\prime}(y)$.
We consider $h$ as a function of $t$.
Since ${\Cal L}$ is in the $c\hvar^{1-\frac{1}{r}}$-cone 
of $\pm\left(\p\sb{x\csb{n}}\right)\sb{\eta\cpr}$,
$$
\p\sb{t}^j h
-\vect{K}\ssb{R}^j h
=O(c\hvar^{1-\frac{1}{r}}).
$$
At the points 
$(x,\vartheta)$ and $(y,\vartheta)$,
$\sigma\sb{j}\vect{K}\ssb{R}\sp{j}h\ge-\hvar$;
we assume $c$ is so small that 
$$
\sigma\sb{j}h\sp{(j)}(t)
\ge -\hvar^{1-\frac{1}{r}}
\qquad\text{at the points } 
t=0 \text{ and } t=\abs{{\Cal L}}.
\tag{\sec5.6}
$$

Since $\abs{\vect{K}\ssb{R}^{r} h}\ge \varkappa\ssb{R}$,
we also know that (again, assuming that $c$ is 
sufficiently small)
$$
\abs{h\sp{(r)}(t)}
\ge \frac{\varkappa\ssb{R}}{2}>0,
\qquad\text{for any } t \text{ between }
0 \text{ and }\abs{{\Cal L}}.
\tag{\sec5.7}
$$

If we assumed that $\hvar=0$
in the right-hand side of $(\sec5.6)$,
then,
similarly to the argument in the proof of Lemma \sec2.1,
we would conclude that
the derivatives $h^{(j)}(t)$ of all orders $j<r$
were monotone functions which did not change 
the signs between $t=0$ and $\abs{{\Cal L}}$.
Moreover, we would derive that
$
\abs{h(\abs{{\Cal L}})-h(0)}
\ge\frac{\varkappa\ssb{R}}{2}\cdot\frac{\abs{{\Cal L}}^{r}}{r!}.
$
Since at the endpoints of ${\Cal L}$
the values of $\sigma\sb{j}h\sp{(j)}$
are only greater than $-\hvar^{1-\frac{1}{r}}$,
there is an error involved;
its magnitude is bounded by 
$O(\hvar^{1-\frac{1}{r}}\abs{{\Cal L}})$.
We arrive at
$$
\abs{h(\abs{{\Cal L}})-h(0)}
\ge\frac{\varkappa\ssb{R}}{2}\cdot\frac{\abs{{\Cal L}}^{r}}{r!}
-\const\hvar^{1-\frac{1}{r}}\abs{{\Cal L}}.
\tag\sec5.8
$$
Therefore, $\abs{{\Cal L}}\le\const\hvar^{\frac{1}{r}}$,
and this proves $(\sec5.5)$.


\Refs\widestnumber\key{1997aaaaa}


\ref\key{Ch\ye{85}}
\by M. Christ
\paper\rm 
Hilbert transforms along curves
\jour\it Ann. of Math. \vol 122 \yr 1985
\pages 575--596
\endref

\ref\key{Co\ye{97}}
\by A. Comech
\paper\rm
Integral operators
with singular canonical relations
\inbook in: \it
Spectral theory, microlocal analysis,
singular manifolds,
\rm M. Demuth, E. Schrohe, 
B.-W. Schulze, and J. Sj\"ostrand, eds
\publ Akademie Verlag \publaddr Berlin \yr 1997
\pages 200--248
\endref

\ref\key{Co\ye{98}}
\by A. Comech
\paper\rm
Optimal estimates
for Fourier integral operators
with one-sided folds
\paperinfo preprint
\yr 1997
\endref

\ref\key{Cu\ye{97}}\by S. Cuccagna
\paper\rm 
$L^2$ estimates for  averaging operators  
along curves with two sided $k$ fold singularities
\jour\it Duke Journal
\vol 89 \yr 1997 \pages 203--216
\endref

\ref\key{GrSe\ye{94}}\by A. Greenleaf and A. Seeger
\paper\rm  Fourier integral operators with fold singularities
\jour\it J. Reine Angew. Math. \vol 455 \yr 1994 
\pages 35--56
\endref

\ref\key{GrSe\ye{97a}}\bysame
\paper\rm  Fourier integral operators with cusp singularities
\paperinfo pre\-print \yr 1997
\endref

\ref\key{GrSe\ye{97b}}\bysame
\paper\rm On oscillatory integral operators
with folding canonical relations
\paperinfo preprint
\yr 1997
\endref

\ref\key{GrU\ye{91}}
\by A. Greenleaf and G. Uhlmann
\paper\rm
Composition of some singular Fourier 
integral operators and estimates for 
restricted X-ray transforms. II
\jour\it  Duke Math. J. \vol  64 \yr 1991 \pages 415--444
\endref

\ref\key{H\"o\ye{85}}\by L. H\"ormander 
\book The analysis of linear 
partial differential operators  IV
\publ Springer-Verlag \yr 1985 
\endref

\ref\key{MeT\ye{85}}\by R.B. Melrose and M.E. Taylor 
\paper\rm Near peak scattering and the corrected
Kirchhoff approximation for a convex obstacle
\jour\it Adv. in Math. \vol 55 \yr 1985 \pages 242--315
\endref 

\ref\key{Mo\ye{65}}
\by B. Morin
\paper\rm 
Canonical forms of the singularities of a differentiable mapping
\jour\it C. R. Acad. Sci. Paris \vol 260 \yr 1965 
\pages 6503--6506
\endref

\ref\key{Pa\ye{91}}\by Y.B. Pan 
\paper\rm Uniform estimates for oscillatory integral operators 
\vol 100\yr 1991\pages 207--220
\jour\it J. Funct. Anal. \endref

\ref\key{PaSo\ye{90}}\by Y.B. Pan and C.D. Sogge
\paper\rm Oscillatory integrals associated to folding 
canonical relations
\jour\it Colloq. Math. \vol 61 \yr1990 \pages 413--419 
\endref

\ref\key{Ph\ye{94}}\by D.H. Phong
\paper\rm Singular integrals and Fourier integral operators
\inbook in: \it 
Essays on Fourier Analysis in honor of Elias M. Stein,
\rm  C. Fefferman, R. Fefferman, and S. Wainger, eds
\publ Princeton Univ. Press\yr 1994 \pages 287--320 
\endref

\ref\key{Ph\ye{95}}\bysame
\paper\rm Regularity of Fourier integral operators
\jour\it Proceedings of the International Congress
of Mathematicians\vol 1, 2\yr 1994 \pages 862--874
\endref

\ref\key{PhSt\ye{91}}\by D.H. Phong and E.M. Stein 
\paper\rm Radon transform and torsion 
\jour\it Internat. Math. Res. Notices\vol 4\yr 1991\pages 49--60
\endref

\ref\key{PhSt\ye{94}}\bysame
\paper\rm Models of degenerate Fourier integral operators 
and Radon transforms 
\jour\it Ann. of Math. \vol 140 \yr 1994 \pages 703--722
\endref

\ref\key{PhSt\ye{97}}\bysame
\paper\rm 
Newton polyhedron and oscillatory integral operators
\paperinfo preprint \yr 1997
\endref

\ref\key{Se\ye{93}}\by A. Seeger
\paper\rm Degenerate Fourier integral operators in the plane
\jour\it Duke Math. J. \vol 71 \yr 1993\pages 685--745
\endref

\ref\key{SoSt\ye{86}}
\by C.D. Sogge and E.M. Stein
\paper\rm  Averages over hypersurfaces II
\jour\it Invent. Math. \vol 86\yr 1986 
\pages 233--242
\endref

\ref\key{St\ye{93}}\by  E.M. Stein
\book Harmonic Analysis: real-variable methods,
orthogonality, and oscillatory integrals 
\publ Prince\-ton Univ. Press \yr 1993 \endref

\noref\key{Th\ye{63}}
\by R. Thom
\paper\rm 
Les singularit\'es
des applications differentiables
\jour\it Ann. Inst. Fourier
\vol 6 \yr 1963 
\pages 43--87
\endref

\endRefs
\enddocument